# FROM CONJECTURE GENERATION BY MAINTAINING DRAGGING TO PROOF

Anna Baccaglini-Frank
Department of Mathematics
"Sapienza" University of Rome (Italy)

Samuele Antonini
Department of Mathematics
University of Pavia (Italy)

*In this paper we propose a hypothesis about how different uses of maintaining dragging, either as a physical tool in a dynamic geometry environment or as a psychological tool for generating conjectures can influence subsequent processes of proving. Through two examples we support the hypothesis that using maintaining dragging as a physical tool may foster cognitive rupture between the conjecturing phase and the proof, while using it as a psychological tool may foster cognitive unity between them.*

## THEORETICAL PERSPECTIVE

Mathematics educators have been encouraging the use of technology in the classroom, and, in particular, several studies on the teaching and learning of geometry (e.g., Noss & Hoyles, 1996; Mariotti, 2006) have shown that a Dynamic Geometry Environment (DGE) can foster the learners' processes of conjecture generation and argumentation, especially in open problem situations, for which the dragging tool plays a crucial role (e.g. Arzarello et al., 2002; Baccaglini-Frank & Mariotti, 2010; Leung, Baccaglini-Frank & Mariotti, 2013). In particular, research carried out by Arzarello et al. (2002) led to the description of some dragging modalities used by secondary school students when asked to solve open problems by producing conjectures in a DGE and then proving them. They describe the key moment of the process of conjecture generation as an *abduction*, related to the use of a particular form of dragging used to maintain a certain geometrical property while the figure changes as an effect of dragging one of its points. We now clarify the main theoretical notions used in the paper.

### Maintaining Dragging

To shed light on the key moment described above, the first author conducted a study (Baccaglini-Frank, 2010a; Baccaglini-Frank & Mariotti, 2010, 2011), in which intentionally inducing an invariant by dragging a point was called *maintaining dragging* (an example will be shown below), and this was explicitly introduced to students between the ages of 15 and 17 in Italian high schools. The study led to a model describing cognitive processes involved in conjecture-generation when maintaining dragging is used by the solver. The findings presented in this paper stem from our interest about possible effects that use of maintaining dragging during conjecture generation might have on the subsequent proof of the conjecture, and were obtained through new analyses of the original data.

## Abduction

Peirce was the first to introduce abduction as the inference, which allows the construction of a claim starting from some data and a rule (Peirce, 1960). In mathematics education there has been renewed interest in the concept of abduction in the context of problem solving in DGEs (e.g., Arzarello et al., 2002; Antonini & Mariotti, 2010; Baccaglini-Frank, 2010b; Baccaglini-Frank & Mariotti, 2011). In this paper we will refer strictly to Piece's definition of abduction of the general form, that is: (*fact*) a fact A is observed; (*rule*) if C were true, then A would certainly be true; (*hypothesis*) so, it is reasonable to assume C is true.

## Cognitive Unity

We will explore relationships between conjecture generation and mathematical proof. Studies have shown that the use of DGEs can promote conjecture generation, but not necessarily the transition to proof (e.g., Yerushalmy, Chazan & Gordon, 1993). However, interesting results have been reached on a possible continuity between these processes, leading to the elaboration of the theoretical construct of *cognitive unity*. The original term (Garuti, Boero & Lemut, 1998) was later redefined, assuming that there may or may not be continuity between the conjecturing phase and the subsequent proof produced (Pedemonte, 2007). The construct of cognitive unity has yielded great potential as a tool of analysis of the relationships between processes of conjecture generation and proofs; cognitive unity can be assessed comparing the sequences of properties logically linked during the conjecturing phase to those elaborated in the proof.

# GENERATING CONJECTURES THROUGH MAINTAINING DRAGGING

We will use an example to show how maintaining dragging can be used to generate conjectures in an open problem situation. The request is the following:

Construct the quadrilateral ABCD (see Fig. 1) following these steps and make conjectures about the possible types of quadrilateral it can become describing *all the ways* you can obtain a particular type of quadrilateral. Construct: a point P and a line *r* through P, the perpendicular line to *r* through P, C on the perpendicular line, a point A symmetric to C with respect to P, a point D on the side of *r* containing A, the circle with centre C and radius CP, point B as the second intersection between the circle and the line through P and D.

The figure can be acted upon by dragging points (let us think about dragging D), and some geometrical properties can be recognized as invariants no matter how the point is dragged (e.g., "CP = PA") while others can become invariants induced through maintaining dragging (e.g., "DA = CB", "CD || BA", "ABCD parallelogram"). During this kind of dragging (e.g., maintaining "ABCD parallelogram" by dragging D), new invariants can be observed as the intentionally induced invariants are visually verified (e.g., "D lies on a circle $C_{AP}$ with centre in A and radius AP").

Recognition of new invariants during dragging can be supported by the use of the trace mark, a functionality in most DGEs. The solver may perceive the newly observed invariants as conditionally linked to the intentionally induced invariant, and express this perception in the form of a conjecture (e.g., "If D belongs to $C_{AP}$, then ABCD is a parallelogram").

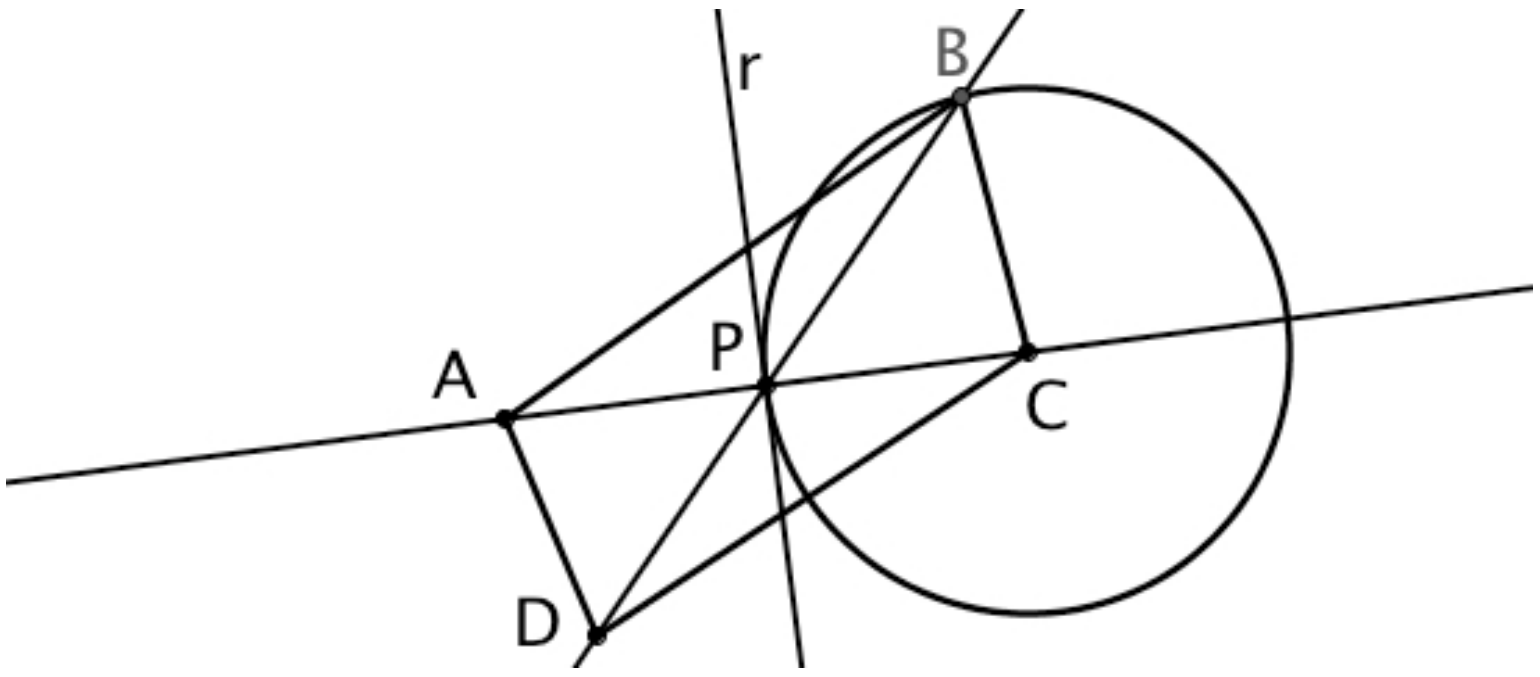

Figure 1: a possible result of the construction in the situation described above.

Our research has shown that many solvers who decide to use maintaining dragging perceive the invariants observed during dragging as *causing* the intentionally induced invariant to be visually verified, at a perceptual level, and interpret this as a *conditional link,* leading to a conjecture in the domain of Euclidean geometry (Baccaglini-Frank & Mariotti, 2010; Leung, Baccaglini-Frank & Mariotti, 2013). The process is described in further detail by Baccaglini-Frank and Mariotti (Baccaglini-Frank & Mariotti, 2011; Baccaglini-Frank, 2010a, 2010b).

## MAINTAINING DRAGGING AS A PSYCHOLOGICAL TOOL

Analyses of the data from the original study have shown that students can also come to use maintaining dragging mentally, freeing it from the physical dragging support. Below is an example of how this happened (also see Baccaglini-Frank 2010a, 2010b from which these excerpts are taken).

Two 15-year-old students in the second year of an Italian high school, Francesco and Gianni, are working on the problem in the example described above. Initially, in order to obtain the desired property (that we will indicate with $P_d$) "ABCD parallelogram" the students have chosen diagonals intersecting at their midpoints ($P_1$) as the property to induce intentionally through maintaining dragging (the student holding the mouse is in bold). However their attempt fails.

Gianni: and now what are we doing? Oh yes, for the parallelogram?

**Francesco**: Yes [as he drags D with the trace activated] yes, we are trying to see when it remains a parallelogram.

Gianni: yes, okay the usual circle comes out.

**Francesco**: wait, because here…oh dear! where is it going? […] So, maybe it's not necessarily the case that D is on a circle so that ABCD is the parallelogram. Because you see, if we then do a kind of circle starting from here, like this, it's good it's good it's good it's good [he drags along a circle he imagines], and then here… see, if I go more or less along a circle that seemed good, instead it's no good…so when is it any good?

Francesco and Gianni give a geometric description of how the point D should be dragged that does not coincide with the trace mark they see on the screen as Francesco performs maintaining dragging. This leads the failure of the students' use of maintaining dragging as a physical tool, so they abandon it. Gianni, who was not dragging, conceives a condition, in his mind. This is shown in the following excerpt.

Gianni: Eh, since this is a chord, it's a chord right? We have to, it means that this has to be an equal cord of another circle, in my opinion with center in A. because I think if you do, like, a circle with center.

**Francesco**: A, you say…

Gianni: symmetric with respect to this one, you have to make it with center A.

[…]

Gianni: with center A and radius AP. I, I think…

**Francesco**: Let's move D. More or less…

Gianni: It looks right doesn't it?

**Francesco**: Yes.

Gianni: Maybe we found it! [Figure 2]

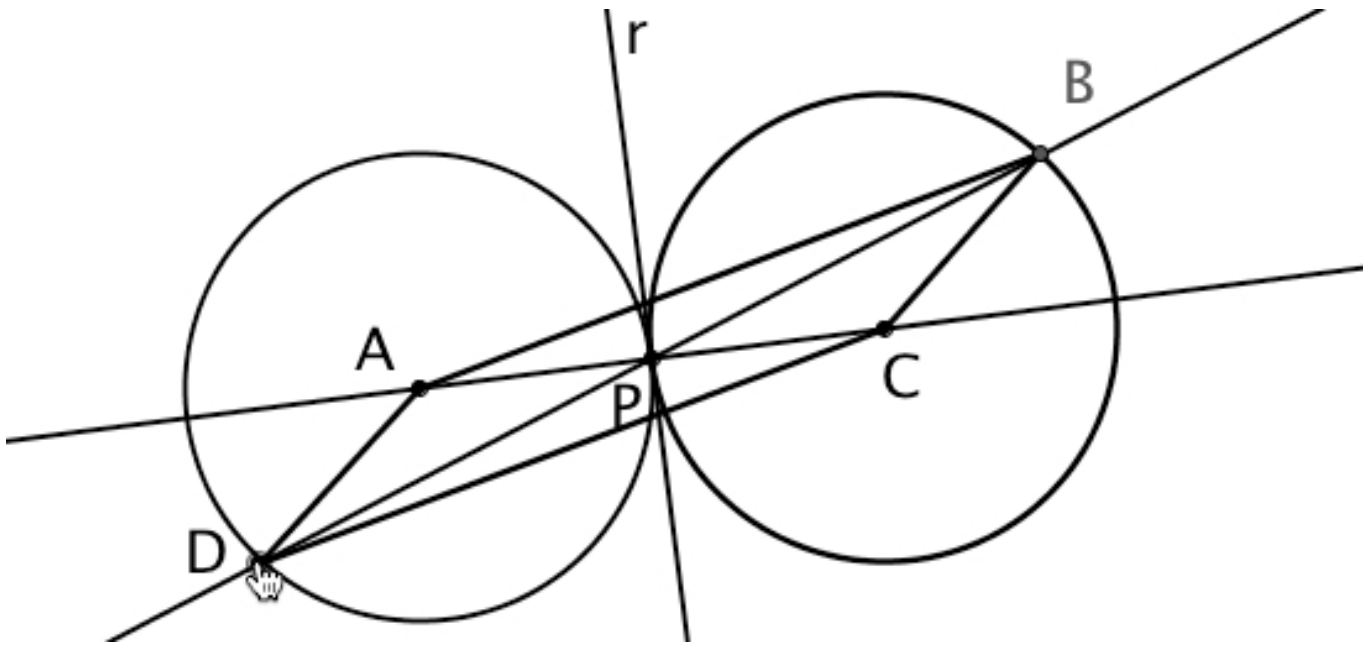

Figure 2: Francesco drags D along the newly constructed circle $C_{AP}$.

Gianni observes that PB is a chord of the circle $C_{CP}$ and reasons abductively:

(*Facts*) PB is the chord of a circle, and PB = PD ($P_2$).

(*Rule*) If PB is the chord of a circle symmetric to the one observed then PB = PD.

(*Abductive hypothesis*) PD is a chord of the symmetric circle $C_{AP}$ ($P_3$).

The abduction leads to a condition as the belonging of D to the circle $C_{AP}$. The students then construct the circle $C_{AP}$ and proceed to link D to it in order to test that when D moves along $C_{AP}$, ABCD is a parallelogram. They seem quite satisfied and formulate the following conjecture immediately after this dragging test: "If D belongs to the circle with center in A and radius AP, then ABCD is a parallelogram".

What happened to maintaining dragging here? The students continue the exploration mentally *as if* they were dragging. Gianni seems to have interiorized the maintaining dragging tool to the point that it as become a *psychological tool* (Vygotsky, 1978, p. 52 ff.) for him (Baccaglini-Frank, 2010a, 2010b). In this case, the conjecturing process relies entirely on his theoretical control over the figure.

## COMPARISON OF TWO PROOFS OF CONJECTURES GENERATED THROUGH DIFFERENT USES OF MAINTAINING DRAGGING

We now compare two proofs of conjectures generated through the two different uses of maintaining dragging described above. Our hypothesis is that when the maintaining dragging tool is physically used, there is a "theoretical gap" left between the premise and the conclusion of the conjecture, that leads to a discontinuity between the conjecturing phase and the proof; on the other hand, if maintaining dragging is used as a psychological tool, the abduction performed by the solver brings out key theoretical ingredients for the proof (similarly to what is described in Arzarello et al., 2002), fostering cognitive unity. The two proofs below are respectively by Gianni and Francesco, the students in the excerpts above, and by Ste and Giu, two 16-year old students.

### Proof in the case of maintaining dragging used as a psychological tool

Below we sketch out the proof constructed by Gianni and Francesco for the conjecture: "If D belongs to the circle with centre in A and radius AP, ABCD is a parallelogram." During the conjecturing phase the following properties were used: ABCD parallelogram ($P_d$); diagonals that intersect at their midpoints ($P_1$); PD = PB ($P_2$); DP chord of a circle symmetric to $C_{CP}$ with centre in A and radius AP ($P_3$); $D \in C_{AP}$ ($P_4$). During this phase the students look at the figure without dragging anything.

The proof was reached in 2 minutes, though the following steps: the circles are symmetric so AD = AP = PC = BC; the isosceles triangles APD and CPB are congruent; so PD = PB; so ABCD has diagonals that intersect at their midpoints, so it is a parallelogram.

| *Conjecturing phase* (time proceeds from right to left) | | | | |
|---|---|---|---|---|
| $P_4$ ← | $P_3$ ← | $P_2$ ← | $P_1$ ← | $P_d$ |
| $D \in C_{AP}$ | DP chord of $C_{AP}$ symmetric to $C_{CP}$ | BP = PD | diagonals intersect at midpoints | ABCD parallelogram |

| *Proof* | | | | | | |
|---|---|---|---|---|---|---|
| $P_4$ and $P_3$ ⇒ | $Q_1$ ⇒ | $Q_2$ ⇒ | $Q_3$ ⇒ | $P_2$ ⇒ | $P_1$ ⇒ | $P_d$ |
| symmetric circles and $D \in C_{AP}$ | AD = AP = PC = BC | ∠DPA= ∠BPC | triangles APD CPB are cong. | BP = PD | diagonals intersect at midpoints | ABCD parallelo gram |

Figure 3: Conjecture and proof generated by Francesco and Gianni.

Figure 3 shows how all the properties used in the conjecturing phase are also used in the proof: the shaded cells show properties used in both phases. In particular, it seems like the abduction generated in the conjecturing phase led to the geometrical properties – $P_1$, $P_2$, $P_3$ – that allowed the students to theoretically fill the gap between the premise ($P_4$) and the conclusion ($P_d$), flipping in particular the abduction between $P_2$ and $P_3$ into a deduction. A few properties – $Q_1$, $Q_2$, $Q_3$ – are added in the proof only in order to theoretically establish the logical relationships between the properties $P_4$, $P_3$ and $P_2$, used in the conjecturing phase. This is a case of cognitive unity between the conjecturing phase and the proof.

**Proof in the case of maintaining dragging used physically**

In order to shed light on relationships between the conjecturing phase and the proof for the second pair of students, we briefly describe how they reached a conjecture. Unlike Gianni and Francesco, Ste and Giu use maintaining dragging physically: they drag D and try to maintain the property "ABCD parallelogram" ($P_d$). Since they have trouble maintaining $P_d$, they induce a different property, that is $B \in C_{PD}$ ($P_3$). Doing this they eventually reach the conjecture: "ABCD is a parallelogram if PA = AD." The properties they touch on during the conjecturing phase are: a parallelogram has diagonals that intersect at their midpoints ($P_1$), BP = PD ($P_2$) through a first abduction, $B \in C_{PD}$ ($P_3$) through a second abduction, $D \in C_{AP}$ ($P_4$) using maintaining dragging, and PA = AD ($P_5$), another invariant observed during maintaining dragging, which they use as premise in their conjecture.

Ste and Giu conclude their proof in about 5 minutes and leave the figure static on the screen as they reason. The inferences they initially make are: CP = PA by construction ($Q_1$), PA = AD ($P_5$), $\angle CPB = \angle APD$ (part of $Q_2$) because vertically opposite angles. These inferences seem not to take into account properties that had been noticed during the exploration. The students are hesitant. Then Giu continues:

Giu: So this is equal to this [he seems to point to BP and PD]…so you prove that one triangle is a 180° rotation of the other and you prove these are parallel? I don't know. […]

Ste: We need to prove that B belongs to this circle here…

Giu seems to be pointing to BP and PD, stating they are equal ($P_2$), however, he is interrupted by Ste, and in the end this property is not used in the proof. The students try (and in the end succeed) to prove that the triangles CPB and APD are isosceles and congruent ($Q_3$): they claim the triangles have congruent sides of the same lengths ($Q_1$) and equal base angles ($\angle CBP = \angle CPB = \angle APD = \angle ADP$) ($Q_2$). Since $\angle CBP = \angle ADP$ (part of $Q_2$) and DA = CB ($Q_4$), then two opposite sides of ABCD are not only congruent but also parallel ($Q_5$), which proves that ABCD is a parallelogram ($P_d$).

Surprisingly to us, in the final proof the students never use the property BP = PD ($P_2$), which played a key role in the conjecturing phase, leading to a property that they found easier to maintain during dragging and that eventually led them to the property they used as a premise. Giu seems to recall the importance of the property as he points to BP and PD at the beginning of the proving phase, but there seems to be too great a “theoretical gap” between $P_4$ and $P_3$, left by the physical use of maintaining dragging, for the students to be able to theoretically bridge the gap and make use of the properties identified in the conjecturing phase. Figure 4 summarizes the properties used in the conjecturing phase and in the proof by Giu and Ste. The shaded cells show the common properties, which end up being only the premise ($P_5$) and the conclusion ($P_d$) of the conditional statement, now a proved theorem. The lack of cognitive unity between the conjecturing phase and the proof is quite evident.

| *Conjecturing phase* (time proceeds from right to left) | | | | | |
|---|---|---|---|---|---|
| $P_5$ ← | $P_4$ ← | $P_3$ ← | $P_2$ ← | $P_1$ ← | $P_d$ |
| PA = AD | D ∈ $C_{AP}$ | B ∈ $C_{PD}$ | BP = PD | diagonals intersect at midpoints | ABCD parallelogram |
| *Proof* | | | | | |
| $P_5$ & | $Q_1$ ⇒ | $Q_2$ ⇒ | $Q_3$ ⇒ | $Q_4$ & $Q_5$ ⇒ | $P_d$ |
| PA = AD | CP = PA | ∠CBP = CPB = ∠APD = ∠ADP | triangles CPB and DPA are cong. | DA = CB and CB ∥ AD | ABCD parallelogram |

Figure 4: Conjecture and proof generated by Giu and Ste

## CONCLUSIONS

Although we cannot draw any general conclusions because we still have very few data on processes of proving for this kind of open problem situation, the two examples offered seem to well support our hypothesis that different uses of maintaining dragging may foster cognitive unity or rupture between the conjecturing phase and the proof of a conjecture generated by students in open problem situations. We find particularly interesting that when maintaining dragging is internalized and becomes a psychological tool, and it is no longer used physically, the abductive reasoning that takes place in the conjecturing phase seems to lead to the discovery of geometrical elements and properties that otherwise are not noticed. These can be reinvested in the proving phase, since they can be re-elaborated into the deductive steps of a proof, as in the case of Gianni and Francesco. On the other hand, in the case of physical use of maintaining dragging these geometrical elements are “absorbed” by the tool: the conjecturing phase seems to not allow the solvers to “bridge the gap” between the premise and the conclusion of their conjecture.

We believe that considerations emerging from this paper can help educators establish educational goals when geometry is taught with the support of DGEs, or at least provide them with food for thought. Important issues to further investigate are whether we want students to become proficient enough in the use of maintaining dragging for it to become a psychological tool for them. If so, we could explore how this might be accomplished in educational settings.